\documentclass[12pt,a4paper]{article}
\usepackage[a4paper,left=2.5cm,right=2.5cm]{geometry}
\usepackage{amsmath,amssymb,amsfonts,amsthm,enumerate,float}
\usepackage{mathrsfs,amsbsy,bm,mathtools,color}
\usepackage{indent first}
\usepackage{scalerel}

\newtheorem{Thm}{Theorem}[section]%
\newtheorem{Lem}[Thm]{Lemma}%
\newtheorem{Prop}[Thm]{Proposition}%
\newtheorem{Cla}[Thm]{Claim}

\theoremstyle{definition} 
\newtheorem{Rmk}[Thm]{Remark}%
\newtheorem{Quest}[Thm]{Question}%

\newcommand{\cN}{\mathcal{N}}

\numberwithin{equation}{section}

\title{Frequently Hypercyclic Meromorphic Curves with Slow Growth}

\author{
  Zhangchi Chen \and
  Bin Guo \and
  Song-Yan Xie
}

\begin{document}
\maketitle

\begin{abstract}
We construct entire curves in projective spaces that exhibit frequent hypercyclicity under translations along countably many prescribed directions while maintaining optimal slow growth rates. Furthermore, we establish a fundamental dichotomy by proving the impossibility of such curves simultaneously preserving frequent hypercyclicity for uncountably many directions under equivalent growth constraints. This result reveals a striking contrast with classical hypercyclicity phenomena, where entire functions can achieve hypercyclicity over some uncountable direction set without growth rate compromise. Our methodology is rooted in Nevanlinna theory and guided by the Oka principle, offering new insights into the relationship between dynamical properties and growth rates of entire curves in projective spaces.
\end{abstract}


\section{Introduction}

In 1929, Birkhoff~\cite{Birkhoff-1929} discovered a striking phenomenon on the space \(\mathcal{H}(\mathbb{C})\) of entire functions endowed with the compact–open topology: for every nonzero \(a \in \mathbb{C}\), the translation operator  
\[
{\sf T}_a \colon \mathcal{H}(\mathbb{C}) \to \mathcal{H}(\mathbb{C}), \qquad
{\sf T}_a f\,(z) = f(z + a),
\]
admits elements that Birkhoff called universal—later termed hypercyclic.  
That is, there exists an entire function \(f\) whose orbit \(\{{\sf T}_a^n f\}_{n \in \mathbb{N}}\) under iteration of \({\sf T}_a\) is dense in \(\mathcal{H}(\mathbb{C})\).  
In modern terminology, such an \(f\) is called hypercyclic for \({\sf T}_a\), and the operator \({\sf T}_a\) itself is said to be hypercyclic.

This discovery revealed that even simple linear maps—such as translations—can generate rich, chaotic dynamics in infinite-dimensional settings, thereby laying the foundation for what is now known as linear chaos or linear dynamical systems.  
For comprehensive treatments of this theory we refer to the monographs~\cite{Bayart-Matheron-2009, Grosse-Erdmann-Manguillot-2011}.

In the present paper, we study hypercyclicity beyond the linear framework.  
Let \(X\) be a topological space and \( {\sf T}\colon X \to X\) a continuous (not necessarily linear) operator.  
An element \(x \in X\) is called hypercyclic for \({\sf T}\) if its orbit \(\{{\sf T}^n x\}_{n \in \mathbb{N}}\) is dense in \(X\).

Let \(\mathcal{H}(\mathbb{C}, \mathbb{P}^m)\) denote the space of holomorphic maps from \(\mathbb{C}\) to the complex projective \(m\)-space \(\mathbb{P}^m\), endowed with the compact-open topology. For a fixed nonzero number \(a \in \mathbb{C}\), the translation operator \({\sf T}_a\) on \(\mathcal{H}(\mathbb{C}, \mathbb{P}^m)\) is defined by \({\sf T}_a(f)(z) = f(z + a)\). It is evident that any hypercyclic entire or meromorphic function \(h \in \mathcal{H}(\mathbb{C}, \mathbb{P}^m)\) (with respect to \({\sf T}_a\)) must be transcendental. In Nevanlinna theory, the complexity of an entire curve \(h: \mathbb{C} \rightarrow \mathbb{P}^m\) is measured using the Nevanlinna–Shimizu–Ahlfors characteristic function (cf., e.g.,~\cite{Noguchi-Winkelmann-2014}):
\begin{equation}\label{eqn-char}
T_h(r) := \int_{t=0}^r \frac{\mathrm{d}t}{t} \int_{|z|<t} h^* \omega_{\sf FS},
\qquad
\forall r \geqslant 0,
\end{equation}
where \(h^* \omega_{\sf FS}\) denotes the pullback by \(h\) of the Fubini–Study form \(\omega_{\sf FS}\) on \(\mathbb{P}^m\). 

Dinh and Sibony posed the following challenge~\cite[Problem~9.1]{Dinh-Sibony-2020}: determine the minimal possible growth rate \(T_h(\cdot)\) for a meromorphic function \(h\colon \mathbb{C}\to \mathbb{P}^1\) that is hypercyclic with respect to some translation operator \({\sf T}_a\). 

It is noted that for any such hypercyclic \(h\), the characteristic function must satisfy

\begin{equation}\label{bigologr}
\lim_{r\to +\infty} \frac{T_h(r)}{\log r} = +\infty .
\end{equation}
Indeed, for a fixed nonconstant target $g \in \mathcal{H}(\mathbb{C}, \mathbb{P}^1)$, there exists a sequence $\{a_n\}$ tending to infinity such that the translates $\mathsf{T}_{a_n}(h)$ converge uniformly to $g$ on a neighborhood of the closed unit disc $\overline{\mathbb{D}}$. Consequently, the image $h(\mathbb{C})$ (counting multiplicities) must contain infinitely many holomorphic discs – obtained from pairwise disjoint discs in the source $\mathbb{C}$ – whose images approximate the fixed holomorphic disc $g(\overline{\mathbb{D}})$. 
Because $\int_{\overline{\mathbb{D}}} g^* \omega_{\sf FS} > 0$, and because the uniform convergence on neighborhoods of the discs ensures that the integral over each translated disc approaches this positive value as $n$ increases, each such disc contributes an area that is asymptotically close to $\int_{\overline{\mathbb{D}}} g^* \omega_{\sf FS}$. Hence the total area 
$
\int_{|z|<t} h^* \omega_{\textsf{FS}}
$
necessarily diverges as $t \to +\infty$; this geometric-topological constraint forces the growth lower bound \eqref{bigologr}.

The following theorem can be seen as reflecting the {Oka principle} — a guiding insight in several complex variables that analytic problems often have only topological obstructions~\cite[pp.~368–369]{Forstneric-2023}.
 In our setting, this principle leads one to expect that the  growth estimate~\eqref{bigologr} might in fact be the sole obstruction to the Dinh–Sibony question.  
Indeed, the theorem below confirms that once \eqref{bigologr} is satisfied, no further analytic barrier remains.

\begin{Thm}[\cite{Chen-Huynh-Xie-2023}--Theorem 1.3]\label{CHX theorem}
For any \( m \geqslant 1 \), for any continuous increasing function \(\phi: [1, +\infty) \to (0, +\infty)\) tending to infinity, and for any countable set \( E \subset [0, 2\pi) \), there exists an entire curve \( h: \mathbb{C} \rightarrow \mathbb{P}^m \) with slow growth rate
\[
T_h(r) \leqslant \phi(r) \cdot \log r,
\qquad
\forall r > 2,
\]
such that \( h \) is simultaneously hypercyclic with respect to every translation \(\mathsf{T}_a\), where \( a \in \mathbb{C} \setminus \{0\} \) has an argument in \( E \).
\end{Thm}

While Theorem~\ref{CHX theorem} establishes the existence of hypercyclic entire curves with optimally slow growth, our investigation advances this program by addressing a refined dynamical question: what minimal growth rates permit the \emph{frequent hypercyclicity} of  entire curves in projective spaces? To articulate this strengthened dynamical behavior, we first ground the discussion with precise definitions.

The \emph{lower density} of a subset \( A \subset \mathbb{N} \), quantifying its asymptotic frequency, is formally defined by
\begin{equation}
    \label{define-lower-density}
  \underline{\mathrm{dens}}(A) \coloneqq
  \liminf_{N \to +\infty} \frac{\# \{ n \leqslant N : n \in A \}}{N}.
\end{equation}

An important development in hypercyclic operator theory was achieved by Bayart and Grivaux~\cite{Bayart-Grivaux-2006}, who introduced the concept of frequent hypercyclicity. This refinement requires a continuous operator \(\mathsf{T}: V \rightarrow V\) to admit elements \( x \in V \) (termed \emph{frequently hypercyclic elements}) such that for every non-empty open subset \( U \subset V \), the visitation set
\[
\{ n \in \mathbb{N} : \mathsf{T}^n(x) \in U \}
\]
has strictly positive lower density. The family of such elements is denoted by \(\mathsf{FHC}(\mathsf{T})\).

The concept of frequent hypercyclicity can be traced back to Voronin's study~\cite{Voronin} of the Riemann zeta function \(\zeta\). Voronin discovered that the translations \(\mathsf{T}_{\sqrt{-1}t} \zeta(\cdot) := \zeta(\cdot + \sqrt{-1}t)\) (for \( t > 0 \)) of the zeta function along the imaginary axis exhibit a remarkable ``universality" property, known as Voronin's theorem. Assuming the Riemann hypothesis, all such translations \(\mathsf{T}_{\sqrt{-1}t}\) of the zeta function \(\zeta\) demonstrate a kind of astonishing frequent hypercyclicity~\cite[p.~265]{Bayart-Matheron-2009}. This result underscores the importance of the concept of frequent hypercyclicity in both complex analysis and operator theory, as it provides a deeper understanding of the dynamical properties of the zeta function and its translations.

The ubiquity of frequently hypercyclic operators and elements in complex geometry is demonstrated by the following result. Recall that a connected complex manifold \( Y \) is called an \emph{Oka manifold} if it satisfies the Oka principle, meaning that holomorphic maps from Stein manifolds to \( Y \) enjoy strong approximation and extension properties; see \cite{Forstneric-2023} for a comprehensive treatment. In particular, both the complex Euclidean space \( \mathbb{C}^m \) and the complex projective space \( \mathbb{P}^m \) are Oka manifolds. 

Let \( Y \) be a connected Oka manifold and let \( n \geqslant 1 \) be an integer. On the space \( \mathcal{H}(\mathbb{C}^n, Y) \) of holomorphic maps from \( \mathbb{C}^n \) to \( Y \), endowed with the compact-open topology, we define for each nonzero vector \( a \in \mathbb{C}^n \) the \emph{translation operator}
\begin{equation}\label{eq:translation-operator}
{\sf T}_a(f)\,(z) = f(z + a), \qquad \forall f \in \mathcal{H}(\mathbb{C}^n, Y), \ \forall z \in \mathbb{C}^n.    
\end{equation}

\begin{Thm}[\cite{Guo-Xie-2}--Theorem D]
Let \( Y \) be a connected Oka manifold, and let \( n \geqslant 1 \) be an integer. Then the set
\[
\bigcap_{b \in \mathbb{C}^n \setminus \{\mathbf{0}\}} \mathsf{FHC}({\sf T}_b)
\]
is dense in \( \mathcal{H}(\mathbb{C}^n, Y) \).
\end{Thm}

\smallskip

The \emph{order} of a holomorphic map \( h \in \mathcal{H}(\mathbb{C}^n, \mathbb{C}^m) \) or \( h \in \mathcal{H}(\mathbb{C}^n, \mathbb{P}^m) \) is defined as
\[
\rho_h := \limsup_{r \to +\infty} \frac{\log T_h(r)}{\log r} \in [0, +\infty],
\]
which quantifies the asymptotic growth rate of \( h \) (see, e.g.,~\cite[(2.3.30)]{Noguchi-Winkelmann-2014} for a definition of \( T_h(r) \) when \( n > 1 \)). For any nonzero vector \( a \in \mathbb{C}^n \), if an entire map \( h \) is frequently hypercyclic for \( \mathsf{T}_a \) (defined in~\eqref{eq:translation-operator}), then a standard application of the First Main Theorem in Nevanlinna theory (cf.~\cite{Noguchi-Winkelmann-2014, Ru-2021}) yields that its order satisfies \( \rho_h \geqslant 1 \). Equivalently, there exists a constant \( C_h > 0 \) such that
\begin{equation}
    \label{optimal hope}
T_h(r) \geqslant C_h \cdot r,
\qquad
\forall r > 1.
\end{equation}
This inequality constitutes a \emph{topological obstruction} on the growth of frequently hypercyclic entire maps. In the spirit of the Oka principle, which posits that topological obstructions are often the only ones in holomorphic settings, it is natural to ask whether the lower bound~\eqref{optimal hope} is essentially optimal, i.e., whether there are no further \emph{analytic} obstructions beyond this necessary condition.

\begin{Quest}[see also~\cite{Guo-Xie-2}--Question 1.14]\label{quest-fhc}
Let \( m, n \geqslant 1 \) be integers and let \( \{a_i\}_{i \in \mathbb{N}} \subset \mathbb{C}^n \setminus \{\mathbf{0}\} \) be a countable set of translation vectors. Consider the translation operators \( \mathsf{T}_{a_i} \) defined as in~\eqref{eq:translation-operator}, acting either on the space \( \mathcal{H}(\mathbb{C}^n, \mathbb{C}^m) \) or on \( \mathcal{H}(\mathbb{C}^n, \mathbb{P}^m) \), both endowed with the compact-open topology. Does there exist a holomorphic map \( h \) in the respective space with \emph{minimal order} \( \rho_h = 1 \), which is simultaneously \emph{frequently hypercyclic} for all operators \( \mathsf{T}_{a_i} \)?
\end{Quest}

\begin{Rmk}
Blasco, Bonilla, and Grosse-Erdmann~\cite[Theorem 3.1]{BBGE} proved the existence of entire functions of one variable in \( \mathsf{FHC}(\mathsf{T}_a) \) of order $1$ for any given \( a \in \mathbb{C} \setminus \{0\} \), which is optimal with respect to the lower bound~\eqref{optimal hope}. 

Their method can be extended to higher dimensions directly: for any integers \( m, n \geqslant 1 \) and any nonzero vector \( a \in \mathbb{C}^n \setminus \{\mathbf{0}\} \), one can construct an entire map \( h \in \mathsf{FHC}(\mathsf{T}_a) \subset \mathcal{H}(\mathbb{C}^n, \mathbb{C}^m) \) of order $1$.

When \( n \leqslant m \), by further applying Sard's theorem as in~\cite[Section~4.3]{Guo-Xie}, one can also obtain elements in \( \mathsf{FHC}(\mathsf{T}_a) \subset \mathcal{H}(\mathbb{C}^n, \mathbb{P}^m) \) of order 1. 

However, the construction in~\cite[Theorem 3.1]{BBGE} — and its extension described above — inherently requires that all translation directions be real multiples of a fixed direction. Consequently, it does not apply directly to the setting of Question~\ref{quest-fhc}, where one aims for simultaneous frequent hypercyclicity for a countable family of translation operators \( \{\mathsf{T}_{a_i}\}_{i \in \mathbb{N}} \) whose directions \( a_i \in \mathbb{C}^n \setminus \{\mathbf{0}\} \) are not necessarily colinear over \( \mathbb{R} \).
\end{Rmk}

The following main result of this article answers Question~\ref{quest-fhc} affirmatively for \( \mathcal{H}(\mathbb{C}, \mathbb{P}^m) \) for any \( m \geqslant 1 \), and solves part of~\cite[Question 1.14]{Guo-Xie-2}. Some inspiration comes from the Oka principle.

\begin{Thm}\label{main-thm-c}
For any integer \( m \geqslant 1 \), for any \( \epsilon > 0 \), and for any countable set \( E \subset [0, 2\pi) \), there exists an entire curve \( h: \mathbb{C} \rightarrow \mathbb{P}^m \) with slow growth rate
\[
T_h(r) \leqslant \epsilon \cdot r,
\qquad
\forall r > 0,
\]
which is simultaneously frequently hypercyclic with respect to every translation \( \mathsf{T}_{R \cdot e^{\sqrt{-1} \cdot \theta}} \), where \( R > 0 \) and \( \theta \in E \).
\end{Thm}

\medskip

 It is crucial to emphasize that the requirement for \( E \) to be a countable set is essential and cannot be relaxed, because of the following:

\begin{Prop}
\label{countable fhc direction}
Let \( h: \mathbb{C} \rightarrow \mathbb{P}^m \) be an entire curve with  slow growth rate
\[
T_h(r) \leqslant \epsilon \cdot r,
\qquad
\forall r > 0,
\]
for some constant \( \epsilon > 0 \).
Then the set $E_h$ of directions \(\theta \in [0, 2\pi)\) for which $h\in{\sf FHC}(\mathsf{T}_{e^{\sqrt{-1}\cdot\theta}})$ is at most countable.
\end{Prop}

\medskip

The complete argument will be developed in Section~\ref{sect: Nevanlinna Theory}. This cardinality restriction reveals a striking dichotomy: while hypercyclic entire functions can maintain minimal growth rates across uncountably many directions~\cite[Theorem D]{Guo-Xie}, our result establishes that frequent hypercyclicity for  entire curves
in projective spaces 
under slow growth conditions necessarily concentrates along countable direction sets. 

\medskip

Our proof strategy for Theorem~\ref{main-thm-c} builds upon the foundational construction in Theorem~\ref{CHX theorem}~\cite{Chen-Huynh-Xie-2023}, which we first outline for the prototypical case \( m = 1 \).

The core construction begins with a carefully curated sequence \(\{\gamma^{[k]}(z)\}_{k \geqslant 1}\) of rational functions vanishing at \( z = \infty \). The hypercyclic meromorphic function \( h \) in Theorem~\ref{CHX theorem} materializes as
\begin{equation}
    \label{magic shape}
h(z) = \sum_{k \in \mathbb{N}} \gamma^{[k]}(z + c_k),
\end{equation}
where the translation parameters \(\{c_k\}_{k \in \mathbb{N}}\) escape to \(\infty\) with controlled rapidity:
\[
1 \ll |c_1| \ll |c_2| \ll |c_3| \ll \cdots.
\]
Through Nevanlinna's First Main Theorem (see Section~\ref{sect: Nevanlinna Theory}), this configuration ensures convergence to a slowly growing hypercyclic meromorphic function, as explained in~\cite[Section~2.3]{Guo-Xie}.

\medskip

Implementing this framework for Theorem~\ref{main-thm-c} introduces a pivotal complication: the sequence \(\{c_n\}_{n \in \mathbb{N}}\) must simultaneously maintain positive lower density along every radial direction \(\mathbb{R}_+ \cdot e^{\sqrt{-1} \cdot \theta}\) for \( \theta \in E \). We address this density preservation challenge through Lemma~\ref{key lem}
(cf.~\cite[Lemma 2.2]{Bayart-Grivaux-2006},~\cite[Lemma 9.5]{Grosse-Erdmann-Manguillot-2011}),
providing a new geometric proof that offers fresh insights into the fundamental mechanism.
Our self-contained approach not only presents an alternative to existing arguments
but also demonstrates potential applicability to a wider range of density-sensitive
constructions (Remark~\ref{potential applicatoin}).

\medskip

Convergence control for \( h \) in~\eqref{magic shape} under the density constraints demands precise calibration between the vanishing orders of \(\gamma^{[k]}(z)\) at infinity and the spatial distribution of \(\{c_n\}_{n \in \mathbb{N}}\). Section~\ref{sect: FHC curves with slow growth} details our solution through  enhancement of vanishing orders combined with careful parameter selection.

\medskip

Finally, establishing the growth bound \( T_h(r) \leqslant \epsilon \cdot r \) employs a scaling transformation in~\eqref{magic trick}. This innovation circumvents intricate Fubini-Study metric computations present in~\cite{Chen-Huynh-Xie-2023, Guo-Xie}, though substantial analytic estimates remain unavoidable (Section~\ref{section 4}).

\section{Nevanlinna Theory}
\label{sect: Nevanlinna Theory}

Let \([Z_0 : \dots : Z_m]\) be the standard homogeneous coordinates of \(\mathbb{P}^m\). Let \(h: \mathbb{C} \rightarrow \mathbb{P}^m\) be an entire curve whose image is not contained in the coordinate hyperplane \(H_0 := \{Z_0 = 0\}\). Fix a reduced representation \(h = [h_0 : h_1 : \dots : h_m]\), i.e., \(h_0, h_1, \dots, h_m\) are holomorphic functions on \(\mathbb{C}\) without any common zero. For simplicity, assume that \(h(0) \notin H_0\).

The \emph{counting function} \(N_h(r, H_0)\) of \(h\) with respect to \(H_0\) is defined as
\[
N_h(r, H_0) := \int_0^r n_h(t, H_0) \, \frac{\mathrm{d}t}{t},
\qquad
\forall r \geqslant 0,
\]
where \(n_h(t, H_0)\) is the number of zeros of \(h_0\) on the disc \(\{|z| < t\}\) counting multiplicities.

The \emph{proximity function} \(m_h(r, H_0)\) of \(h\) with respect to \(H_0\) is given by
\[
m_h(r, H_0) := \int_{\theta=0}^{2\pi} \frac{1}{2} \log \left(1 + \left|\frac{h_1}{h_0}\right|^2 (re^{\sqrt{-1} \cdot \theta}) + \dots + \left|\frac{h_m}{h_0}\right|^2 (re^{\sqrt{-1} \cdot \theta})\right) \, \frac{\mathrm{d}\theta}{2\pi},
\qquad
\forall r \geqslant 0.
\]

One version of the First Main Theorem in Nevanlinna theory (cf.~\cite[Theorem 2.3.31]{Noguchi-Winkelmann-2014}) states that the Nevanlinna-Shimizu-Ahlfors characteristic function~\eqref{eqn-char} can be rewritten as
\begin{align}\label{fmt}
T_h(r) = m_h(r, H_0) + N_h(r, H_0) - m_h(0, H_0),
\qquad
\forall r \geqslant 0.
\end{align}

\begin{Prop}\label{optimal slow rate}
If an entire curve \(h: \mathbb{C} \rightarrow \mathbb{P}^m\) is frequently hypercyclic for some translation \(\mathsf{T}_a\), where \(a \in \mathbb{C} \setminus \{0\}\), then there exist constants \(C > 0\) and \(R_0 > 1\) such that
\[
T_h(r) \geqslant C \cdot r,
\qquad
\forall r \geqslant R_0.
\]
\end{Prop}

The proof relies on a classical principle of complex analysis, already present in the proof of the Hadamard factorization theorem: the growth of an entire function controls the number of its zeros in a disc of radius $r$. In our setting, the frequent hypercyclicity implies that the translates $h(z + n_k a)$ converge uniformly on the unit disc to a fixed nonconstant holomorphic disc $g(\mathbb{D})$ intersecting a hyperplane $H_0 \subset \mathbb{P}^m$; consequently $h$ itself must intersect $H_0$ in infinitely many disjoint discs. Hence the counting function satisfies $N_h(r,H_0) \geqslant C\,r$ for large $r$ for some constant $C>0$, and the First Main Theorem $T_h(r) \geqslant N_h(r,H_0)$ yields the required linear lower bound for $T_h(r)$.

For entire functions $\mathbb{C} \to \mathbb{C}$, the same idea—expressed via Jensen’s formula—produces the growth restriction stated in \cite[Theorem 3.1(b)]{BBGE}. Our proposition reformulates this classical principle in the framework of Nevanlinna theory for entire curves $\mathbb{C} \to \mathbb{P}^m$; we include the short proof for completeness.

\begin{proof}
Consider an entire curve \(g(z) := [z : 1 : \dots : 1] \in \mathcal{H}(\mathbb{C}, \mathbb{P}^m)\) with \(n_g(|a|/3, H_0) = 1\). By Rouché's theorem, there exists a small neighborhood \(\mathcal{U}\) of \(g\) such that \(n_f(|a|/3, H_0) = 1\) for any \(f \in \mathcal{U}\).

Define \(A := \{n \in \mathbb{N} : \mathsf{T}_a^n(h) \in \mathcal{U}\}\). Since \(h \in \mathsf{FHC}(\mathsf{T}_a)\), we have
\[
\alpha := \underline{\mathrm{dens}}(A) = \liminf_{N \to +\infty} \frac{\# \{n \leqslant N : n \in A\}}{N} > 0.
\]
Choose a large \(N_0 \gg 1\) such that
\[
\# \{n \leqslant N : n \in A\} \geqslant \frac{\alpha}{2} \cdot N,
\qquad
\forall N \geqslant N_0.
\]
The disc \(\mathbb{D}(0, (N+1)|a|)\) contains disjoint small discs \(\mathbb{D}(na, |a|/3)\) for \(n \leqslant N, n \in A\). Each such disc contributes exactly one ``zero'' to \(n_h((N+1)|a|, H_0)\). Therefore,
\[
n_h((N+1)|a|, H_0) \geqslant \# \{n \leqslant N : n \in A\} \geqslant \frac{\alpha}{2} \cdot N,
\qquad
\forall N \geqslant N_0.
\]
Thus, we can take \(C = \alpha / 8 > 0\), so that
\begin{equation}
    \label{counting linear bound}
n_h(t, H_0) \geqslant 2C \cdot t,
\qquad
\forall t \geqslant (N_0 + 1)|a|.
\end{equation}
Finally, by the First Main Theorem in Nevanlinna theory, for \(r \geqslant 2(N_0 + 1)|a|\), we have
\[
\begin{aligned}
T_h(r) &\geqslant N_h(r, H_0) \\
&\geqslant \int_{t=(N_0+1)|a|}^r n_h(t, H_0) \, \frac{\mathrm{d}t}{t} \\
\text{[use~\eqref{counting linear bound}]}\quad &\geqslant 2C \cdot (r - (N_0 + 1)|a|) \\
&\geqslant C \cdot r.
\end{aligned}
\]
This completes the proof.
\end{proof}

\smallskip

\begin{proof}[Proof of Proposition~\ref{countable fhc direction}]
Using the same notation as in the preceding proof (with \(|a|=1\)), for each direction \(\theta \in E\), define
\[
A_{\theta} := \{n \in \mathbb{N} : \mathsf{T}_{e^{\sqrt{-1}\theta}}^n(h) \in \mathcal{U}\}.
\]
Since \(h\) is frequently hypercyclic with respect to \(\mathsf{T}_{e^{\sqrt{-1}\theta}}\), we have
\[
\alpha_{\theta} := \underline{\mathrm{dens}}(A_{\theta}) > 0.
\]

We claim that, for each integer \(n \geqslant 1\), the cardinality of the set
\[
E_n := \left\{ \theta \in E : \alpha_{\theta} > \frac{1}{n} \right\}
\]
is less than \(4\epsilon n\). Consequently, since \(E = \bigcup_{n \geqslant 1} E_n\), the set \(E\) is countable.

Suppose, for the sake of contradiction, that \(E_n\) contains \(k \geqslant 4\epsilon n\) distinct directions \(\theta_1, \dots, \theta_k\). We can choose a large \(R \gg 1\) such that the distances between \(R \cdot e^{\sqrt{-1}\theta_i}\) and \(R \cdot e^{\sqrt{-1}\theta_j}\) are pairwise greater than 2024 for all \(i \neq j\).

Furthermore, we can choose a large \(N_1 \gg 1\) such that for all \(N \geqslant N_1\) and \(j = 1, \dots, k\),
\[
\# \{ R\leqslant n \leqslant N-1 : n \in A_{\theta_j} \} \geqslant \frac{\alpha_{\theta_j}}{2} \cdot N.
\]

For each such \( n \), the disc \(\mathbb{D}(n \cdot e^{\sqrt{-1}\theta_j}, 1/3)\) contributes exactly one ``zero'' to the count in \( n_h(N, H_0) \). Noting that these discs are disjoint, we have
\[
n_h(N, H_0) \geqslant \sum_{j=1}^k \frac{\alpha_{\theta_j}}{2} \cdot N \geqslant 4\epsilon n \cdot \frac{1}{2n} \cdot N = 2\epsilon N.
\]

Thus,
\[
N_h(r, H_0) \geqslant \int_{t=N_1}^r n_h(t, H_0) \, \frac{\mathrm{d}t}{t} \geqslant \int_{t=N_1}^r 2\epsilon [t] \, \frac{\mathrm{d}t}{t} > \epsilon \cdot r
\]
for all sufficiently large \(r\).

This contradicts the fact that \(N_h(r, H_0) \leqslant T_h(r) \leqslant \epsilon \cdot r\) for all \(r > 0\). Therefore, the set \(E\) must be at most countable.
\end{proof}

\section{Preparations}
\label{sect: FHC curves with slow growth}

\subsection{A Key Lemma}

The following technical lemma will play a pivotal role in our subsequent analysis. We note that a similar result was originally established by Bayart and Grivaux \cite[Lemma 2.2]{Bayart-Grivaux-2006}.

\begin{Lem}[\cite{Grosse-Erdmann-Manguillot-2011}--Lemma 9.5]\label{key lem}
There exists a family of pairwise disjoint subsets $A(l, \nu)$ of $\mathbb{N}$, indexed by $l, \nu \in \mathbb{N}$, satisfying:
\begin{itemize}
  \item Each $A(l, \nu) \subset \mathbb{N}$ has positive lower density.
  \item For any distinct elements $n_1 \neq n_2 \in \bigcup_{l, \nu \in \mathbb{N}} A(l, \nu)$ with $n_i \in A(l_i, \nu_i)$, one has $n_i \geqslant \nu_i$ and
    \begin{equation}\label{key gap}
      |n_1 - n_2| \geqslant \nu_1 + \nu_2.
    \end{equation}
\end{itemize}
\end{Lem}

Several refined versions of Lemma~\ref{key lem} can be found in the literature, including \cite[Theorem~5.15]{BCP-2021}, \cite[Lemma~2.2]{BGMM-2024}, \cite[Lemma 3.2]{EM-2019}, and \cite[Section 2]{EM-2021}. Although Lemma~\ref{key lem} is well-documented in existing works, we provide a novel constructive proof here for completeness. 

Our methodology is divided into two key steps. First, Proposition~\ref{prop 3.2} establishes the fundamental separation property for parameterized families. Second, Proposition~\ref{prop 3.3} enables a density-preserving decomposition. By combining these results, we derive the  structure presented in Lemma~\ref{key lem}. 

The geometric intuition underlying our construction may have broader implications for related problems, as highlighted in Remark~\ref{potential applicatoin}.

\begin{Prop}\label{prop 3.2}
There exists a family of pairwise disjoint subsets $\{A(\nu)\}_{\nu\in\mathbb{N}}$ in $\mathbb{N}$ satisfying:
\begin{itemize}
  \item Each $A(\nu)$ has positive lower density.
  \item For distinct elements $n_i \in A(\nu_i)$ ($i=1,2$),
    \begin{equation}\label{key gap 0}
      n_i \geqslant \nu_i \quad \text{and} \quad |n_1 - n_2| \geqslant \nu_1 + \nu_2.
    \end{equation}
\end{itemize}
\end{Prop}

\begin{proof}
Our construction employs a two-scale geometric approach. First, fix:
\begin{itemize}
\item A family of pairwise disjoint intervals $\{I_j=[a_j,b_j]\}_{j\geqslant 1}$ in $[0,1]$ with $b_j>a_j$.
\item A family of disjoint positive intervals $\{N_k=[c_k,d_k]\}_{k\geqslant 1}$, 
$0<c_k<d_k<c_{k+1}<d_{k+1}$,
satisfying:
  \begin{equation}\label{expanding trick}
  \lim_{k\to\infty}(d_k-c_k) = \lim_{k\to\infty}(c_{k+1}-d_k) = +\infty,
  \end{equation}
  \begin{equation}\label{has postive density}
  \left(\bigcup_{k\geqslant 1}N_k\right)\cap\mathbb{N} \text{ has positive lower density} \ \delta>0.
  \end{equation}
\end{itemize}

We inductively construct thresholds $\{K_j\}_{j\geqslant 1}$ ensuring:
\begin{equation}\label{gap 11}
c_k - d_{k-1} \geqslant 2j, \quad \forall k\geqslant K_j,
\end{equation}
\begin{equation}\label{gap 2}
\mathrm{dist}(I_\ell,I_j)\cdot|N_k| \geqslant 2j, \quad \forall k\geqslant K_j,\ 1\leqslant\ell<j,
\end{equation}
where $\mathrm{dist}(I_\ell, I_{j})$ denotes the Euclidean distance between the two intervals. 
The existence of such $K_j$ follows from the expansion rates in~\eqref{expanding trick}.

Define scaled intervals of integers:
\[
[I_j\cdot N_k] := [c_k + (d_k-c_k)a_j, c_k + (d_k-c_k)b_j] \cap \mathbb{N},
\quad \forall j, k\geqslant 1,
\]
and form preliminary sets:
\[
A'(\nu) := \bigcup_{k\geqslant K_\nu} [I_\nu\cdot N_k],
\quad \forall \nu\geqslant 1.
\]
Each $A'(\nu)$ inherits positive density from~\eqref{has postive density} through the scaling construction.

For $\nu_1 < \nu_2$ and $n_i \in A'(\nu_i)$ ($i=1, 2$), the separation condition follows from:
\[
|n_1 - n_2| \geqslant 2\nu_2 > \nu_1 + \nu_2
\]
established via~\eqref{gap 11} and~\eqref{gap 2}.

To enforce~\eqref{key gap 0} within each individual family $A(\nu)$, we proceed as follows. Enumerate $A'(\nu) \cap \mathbb{N}_{\geqslant \nu}$ increasingly as $a'_{\nu,1} < a'_{\nu,2} < a'_{\nu,3} < \cdots$, then define
\[
A(\nu) := \{a'_{\nu,2\nu\cdot k}\}_{k \geqslant 1}.
\]
This construction guarantees $|n_1 - n_2| \geqslant 2\nu$ for any distinct $n_1, n_2 \in A(\nu)$. 

By~\eqref{has postive density}, it is clear that each $A'(\nu) \cap \mathbb{N}_{\geqslant \nu}\subset \mathbb{N}$ has positive lower density of at least $|I_{\nu}| \delta>0$. Therefore, $A(\nu)$ inherits a lower density of at least $|I_{\nu}| \delta/(2\nu)>0$ by the following Lemma~\ref{inherit density}.
\end{proof}

\begin{Rmk}
\label{potential applicatoin}
A construction of $\{N_k\}_{k\geqslant 1}$ achieving $\delta=1$ is elementary within this framework. For instance, one may take $c_1:=1$ and define recursively:
\[
d_k - c_k := 2^k, \quad c_{k+1} - d_k := k, \quad \forall k \geqslant 1.
\]

Moreover, for any prescribed $\epsilon > 0$, the intervals $\{I_j\}_{j\geqslant 1}$ can be constructed with precise length decay $|I_{\nu}| = O(\nu^{-(1+\epsilon)})$ as $\nu \to \infty$. Specifically, let
\[
C := \sum_{\nu=1}^\infty \frac{1}{\nu^{1+\epsilon}} < +\infty,
\]
and define preliminary intervals:
\[
I'_\nu := \left[\frac{1}{C}\sum_{j=1}^{\nu-1} \frac{1}{j^{1+\epsilon}}, \frac{1}{C}\sum_{j=1}^\nu \frac{1}{j^{1+\epsilon}}\right], \quad \forall \nu \geqslant 1.
\]
These satisfy $\bigcup_{\nu\geqslant 1} I'_\nu = [0,1]$ with pairwise intersections limited to endpoints. The desired disjoint intervals $\{I(\nu)\}_{\nu\geqslant 1}$ can be obtained by scaling each $I'(\nu)$ by a factor of $\frac{1}{2}$ about its center, yielding:
\[
|I(\nu)| = \frac{1}{2C} \cdot \frac{1}{\nu^{1+\epsilon}}.
\] 

Consequently, the constructed sets $A(\nu) \subset \mathbb{N}$ achieve lower densities of order $O(\nu^{-(2+\epsilon)})$. This controlled decay framework may prove instrumental in other contexts requiring sparse arithmetic configurations.
\qed
\end{Rmk}

\begin{Lem}\label{inherit density}
For two strictly increasing sequences  $\{a_n\}_{n\geqslant 1},\,\{n_k\}_{k\geqslant 1}\subset\mathbb{N}$, if $\underline{\mathrm{dens}}\{a_n\}_{n\geqslant 1}\geqslant \delta_1>0$ and $\underline{\mathrm{dens}}\{n_k\}_{k\geqslant 1}\geqslant\delta_2>0$, then $\underline{\mathrm{dens}}\{a_{n_k}\}_{k\geqslant 1}\geqslant \delta_1\delta_2$.
\end{Lem}

\begin{proof}
By \eqref{define-lower-density}, we have:
\[
\liminf_{N \to \infty} \frac{\# \{n : a_n \leqslant N\}}{N} \geqslant \delta_1 \quad \text{and} \quad \liminf_{N \to \infty} \frac{\# \{k : n_k \leqslant N\}}{N} \geqslant \delta_2.
\]

Fix $0 < \epsilon \ll 1$. For $N \gg 1$, define:
\[
K(N) \coloneqq \max\{n : a_n \leqslant N\}, \quad K'(N) \coloneqq \max\{k : n_k \leqslant K(N)\}.
\]
From the lower density hypotheses:
\[
\frac{K(N)}{N} \geqslant \delta_1 - \epsilon>0 \quad \text{and} \quad \frac{K'(N)}{K(N)} \geqslant \delta_2 - \epsilon>0.
\]

Multiplying these inequalities yields:
\[
\frac{K'(N)}{N} \geqslant (\delta_1 - \epsilon)(\delta_2 - \epsilon).
\]

Taking $\liminf_{N \to \infty}$ followed by $\epsilon \to 0^{+}$, we conclude:
\[
\underline{\mathrm{dens}}\{a_{n_k}\} \geqslant \delta_1 \delta_2>0. \qedhere
\]
\end{proof}

\begin{Prop}\label{prop 3.3}
Any subset $A\subset\mathbb{N}$ with positive lower density admits a decomposition
\[
A = \bigcup_{l\in\mathbb{N}} A_l
\]
into countably many pairwise disjoint subsets $A_l$ with positive lower density.
\end{Prop}

\begin{proof}
For $A=\mathbb{N}$, consider the arithmetic decomposition:
\[
B_\ell := p_\ell\mathbb{N} \setminus \bigcup_{j<\ell}p_j\mathbb{N},  \quad\forall\ell\geqslant 1,
\]
where $\{p_j\}_{j\geqslant 1}$ enumerates primes in  increasing order. Each $B_\ell$ contains positive integers with least prime factor $p_\ell$, having positive lower density. The singleton $\{1\}$ is included in $B_1$ for completeness.

For general $A = \{e_1 < e_2 <e_3< \cdots\}$, project onto the arithmetic structure:
\[
A_l := \{e_k : k \in B_l\},
\quad \forall l\geqslant 1.
\]
By Lemma~\ref{inherit density}, the positive lower density of $A_l$ inherits from the original lower density of $A$ and the lower density of the $B_l$.
\end{proof}

\subsection{Some Estimates}

Let $A(k, l)$ be as defined above in Lemma~\ref{key lem}, where $k, l \in \mathbb{N}$. For each $k \in \mathbb{N}$, choose a large integer (recall \eqref{decay} and \eqref{num poles})
\begin{equation}\label{choose eta}
  \eta_k \geqslant \max\left\{2R_k, 3^k, n_k, \left(2^{k+4} \cdot \sum_{n=1}^{+\infty} \frac{1}{n^3} \right)^{\frac{1}{3}}\right\}.
\end{equation}
For each fixed $k \in \mathbb{N}$, enumerate the set $\bigcup_{l \geqslant \eta_k} A(k, 2l)$ in increasing order as $\{a_s^{[k]}\}_{s \geqslant 1}$. For each $s \in \mathbb{N}$, there is a unique $l_s \geqslant \eta_k$ such that $a_s^{[k]} \in A(k, 2l_s)$. We define $\Theta^{[k]}(s) := l_s$ to record this index.

Fix an auxiliary bijection
\[
  \varphi = (\varphi_1, \varphi_2): \mathbb{N} \overset{\sim}{\longrightarrow} \mathbb{N} \times \mathbb{N}.
\]
We decompose the sequence $\{a_s^{[k]}\}_{s \geqslant 1}$ into countably many disjoint subsets
\begin{equation}\label{center1}
  \mathcal{A}_v^{[k]} := \bigcup_{\substack{l \geqslant \eta_k \\ l \in \varphi_1^{-1}(v)}} A(k, 2l) \qquad \text{for all } v \geqslant 1.
\end{equation}

For a subset \(A \subset \mathbb{N}\) and a complex number \(\lambda \in \mathbb{C}\), we define \(A \cdot \lambda:=\{\,a \lambda : a \in A \,\}\). Enumerate the given countable set \( E \subset [0, 2\pi) \) in Theorem~\ref{main-thm-c} as a sequence $\{\theta(v)\}_{v \geqslant 1}$. Rotate $\mathcal{A}_v^{[k]}$ to $\mathcal{B}_v^{[k]} := \mathcal{A}_v^{[k]} \cdot e^{i \theta(v)}$ for each $v \in \mathbb{N}$. Denote
\begin{equation}\label{set center}
  b_s^{[k]} := a_s^{[k]} e^{i \theta(\varphi_1 \circ \Theta^{[k]}(s))} \qquad \text{for all } k, s \geqslant 1.
\end{equation}
Thus, we rewrite
\begin{equation}\label{center2}
  \mathcal{B}_v^{[k]} = \{b_s^{[k]}: \varphi_1(\Theta^{[k]}(s)) = v\} = \{b_s^{[k]}\}_{s \geqslant 1} \cap \mathbb{N} \cdot e^{\sqrt{-1} \cdot \theta(v)} \qquad \text{for all } k, v \geqslant 1.
\end{equation}
Since $|b_s^{[k]}| = a_s^{[k]} \in A(k, 2\Theta^{[k]}(s))$, by Lemma~\ref{key lem}-(2), we have
\begin{equation}\label{farfrom0}
  |b_s^{[k]}| \geqslant 2\Theta^{[k]}(s) \geqslant 2\eta_k.
\end{equation}

Now, we demonstrate that any two distinct points in $\{b_s^{[k]}\}_{k,s \geqslant 1}$ are at a significant distance from each other.

For convenience, define $b_0^{[k]} = a_0^{[k]} \coloneqq 0$ for each $k \in \mathbb{N}$.

\begin{Prop}\label{distance estimate}
\begin{itemize}
\item[(1)] The discs $\{\mathbb{D}(b_s^{[k]}, \eta_k)\}_{s, k \geqslant 1}$ are pairwise disjoint. Moreover, the annuli
\[
\{|b_s^{[k]}| - \eta_k < |z| < |b_s^{[k]}| + \eta_k\} \quad \text{for all } k, s \geqslant 1
\]
are also pairwise disjoint.

\smallskip
\item[(2)] For any $s, k \geqslant 1$ and $l = \Theta^{[k]}(s) \geqslant \eta_k$, one has
\begin{equation}\label{same k ineq}
\mathsf{dist}(b_t^{[k]}, \overline{\mathbb{D}}(b_s^{[k]}, l)) \geqslant l + 2\eta_k |t - s| \quad \text{for all } t \neq s.
\end{equation}
Furthermore, for any $k' \neq k$, there is a unique integer $\tilde{s} \geqslant 1$ such that $|b_{\tilde{s}-1}^{[k']}| < |b_s^{[k]}| < |b_{\tilde{s}}^{[k']}|$, and
\[
\mathsf{dist}(b_t^{[k']}, \overline{\mathbb{D}}(b_s^{[k]}, l)) \geqslant l + \eta_{k'} (|t - \tilde{s}| + 1) \quad \text{for all } t \geqslant 1.
\]

\item[(3)] For any $b \notin \overline{\mathbb{D}}(b_s^{[k]}, \eta_k)$, there is a unique integer $\tilde{t} \geqslant 1$ such that $|b_{\tilde{t}-1}^{[k]}| \leqslant |b| < |b_{\tilde{t}}^{[k]}|$, and
\begin{equation}\label{outsides disc pt}
|b - b_s^{[k]}| \geqslant \frac{\eta_k}{2} (|\tilde{t} - s| + 1).
\end{equation}
\end{itemize}
\end{Prop}

\begin{proof}
(1) For distinct points $b_s^{[k]}$ and $b_{s'}^{[k']}$, we have
\[
\big||b_s^{[k]}| - |b_{s'}^{[k']}|\big| = |a_s^{[k]} - a_{s'}^{[k']}| \geqslant 2\Theta^{[k]}(s) + 2\Theta^{[k']}(s') \geqslant 2\eta_k + 2\eta_{k'},
\]
which implies that the annuli $\{|b_s^{[k]}| - \eta_k < |z| < |b_s^{[k]}| + \eta_k\}$ and $\{|b_{s'}^{[k']}| - \eta_{k'} < |z| < |b_{s'}^{[k']}| + \eta_{k'}\}$ are disjoint.

\smallskip
\noindent
(2) Since $\{a_j^{[k]}\}_{j \geqslant 1}$ is increasing and $\Theta^{[k]}(j) \geqslant \eta_k$ for all $j, k \geqslant 1$, by \eqref{key gap} we get
\begin{equation}\label{gap 1}
a_{j+1}^{[k]} - a_j^{[k]} \geqslant 2\Theta^{[k]}(j+1) + 2\Theta^{[k]}(j) \geqslant 2\max\{\Theta^{[k]}(j+1), \Theta^{[k]}(j)\} + 2\eta_k.
\end{equation}
Hence, for any two positive integers $t < s$, we have
\[
|b_s^{[k]} - b_t^{[k]}| \geqslant |b_s^{[k]}| - |b_t^{[k]}| = a_s^{[k]} - a_t^{[k]} = \sum_{i=t}^{s-1} (a_{i+1}^{[k]} - a_i^{[k]}) \overset{\eqref{gap 1}}{\geqslant} 2\Theta^{[k]}(s) + 2\eta_k (s - t).
\]
Similarly, for any two positive integers $t > s$, we obtain
\[
|b_s^{[k]} - b_t^{[k]}| \geqslant a_t^{[k]} - a_s^{[k]} \geqslant 2\Theta^{[k]}(s) + 2\eta_k (t - s).
\]
Therefore, \eqref{same k ineq} holds since $l = \Theta^{[k]}(s)$ and
\[
\mathsf{dist}(b_t^{[k]}, \overline{\mathbb{D}}(b_s^{[k]}, l)) \geqslant |b_t^{[k]} - b_s^{[k]}| - l \geqslant l + 2\eta_k |t - s|.
\]

Next, for $k' \neq k$, let $\tilde{s} := \min\{t \in \mathbb{N} : a_t^{[k']} > a_s^{[k]}\}$. Then,
\begin{equation}\label{ktok'}
a_{\tilde{s}-1}^{[k']} < a_s^{[k]} < a_{\tilde{s}}^{[k']}.
\end{equation}
For each $t \geqslant 1$, we have
\[
\mathsf{dist}(b_t^{[k']}, \overline{\mathbb{D}}(b_s^{[k]}, l)) \geqslant |b_t^{[k']} - b_s^{[k]}| - l \geqslant |a_t^{[k']} - a_s^{[k]}| - l.
\]
We now estimate $|a_t^{[k']} - a_s^{[k]}|$.

If $a_t^{[k']} > a_s^{[k]}$, then $a_t^{[k']} \geqslant a_{\tilde{s}}^{[k']}$ and $t \geqslant \tilde{s}$, due to \eqref{ktok'} and the strict increase of $\{a_j^{[k']}\}_{j \geqslant 1}$. Hence,

\begin{align*}
|a_t^{[k']}-a_s^{[k]}|
&=
(a_t^{[k']}-a^{[k']}_{\tilde{s}})+(a^{[k']}_{\tilde{s}}-a^{[k]}_s)
\\
&
=\sum_{i=\tilde{s}}^{t-1}(a^{[k']}_{i+1}-a^{[k']}_{i})+(a^{[k']}_{\tilde{s}}-a^{[k]}_s)\\
[\text{by}~\eqref{key gap}]\qquad&\geqslant\sum_{i=\tilde{s}}^{t-1}\big(2\Theta^{[k']}(i+1)+2\Theta^{[k']}(i)\big)+(2\Theta^{[k]}(s)+2\Theta^{[k']}(\tilde{s}))\\
[\text{remind }\Theta^{[k']}(m)\geqslant\eta_{k'}]\qquad&\geqslant
4\eta_{k'}(t-\tilde{s})+(2l+2\eta_{k'})\\
&\geqslant
\eta_{k'}(|\tilde{s}-t|+1)+2l.
\end{align*}

If $a_t^{[k']} < a_s^{[k]}$, since $\{a_j^{[k']}\}_{j \geqslant 1}$ is strictly increasing, by \eqref{ktok'} we have $a_t^{[k']} \leqslant a_{\tilde{s}-1}^{[k']}$ and $1 \leqslant t \leqslant \tilde{s} - 1$. Therefore,

\begin{align*}
|a_t^{[k']}-a_s^{[k]}|
&
=
(a^{[k']}_{\tilde{s}-1}-a_t^{[k']})+(a^{[k]}_s-a^{[k']}_{\tilde{s}-1})
\\
&=\sum_{i=t}^{\tilde{s}-2}(a^{[k']}_{i+1}-a^{[k']}_{i})+(a^{[k]}_s-a^{[k']}_{\tilde{s}-1})\\
[\text{by }\eqref{key gap}]\qquad
&\geqslant
\sum_{i=t}^{\tilde{s}-2}\big(2\Theta^{[k']}(i+1)+2\Theta^{[k']}(i)\big)+(2\Theta^{[k]}(s)+2\Theta^{[k']}(\tilde{s}-1))\\
[\text{since }\Theta^{[k']}(m)\geqslant\eta_{k'}]\qquad
&\geqslant
4\eta_{k'}(\tilde{s}-1-t)+(2l+2\eta_{k'})\\
&\geqslant\eta_{k'}(|\tilde{s}-t|+1)+2l.
\end{align*}

Summarizing, we conclude that
\[
\mathsf{dist}(b_t^{[k']}, \overline{\mathbb{D}}(b_s^{[k]}, l)) \geqslant |a_t^{[k']} - a_s^{[k]}| - l \geqslant l + \eta_{k'} (|\tilde{s} - t| + 1) \quad \text{for all } t \geqslant 1.
\]

\smallskip
\noindent
(3) Select the unique integer $\tilde{t} \geqslant 1$ such that $a_{\tilde{t}-1}^{[k]} \leqslant |b| < a_{\tilde{t}}^{[k]}$.

If $a_s^{[k]} < |b|$, then $a_s^{[k]} \leqslant a_{\tilde{t}-1}^{[k]}$ and $s \leqslant \tilde{t} - 1$. Hence,
\[
|b - b_s^{[k]}| \geqslant |b| - a_s^{[k]} = (|b| - a_{\tilde{t}-1}^{[k]}) + (a_{\tilde{t}-1}^{[k]} - a_s^{[k]}) \geqslant a_{\tilde{t}-1}^{[k]} - a_s^{[k]} \geqslant 4\eta_k (\tilde{t} - 1 - s).
\]

If $|b| \leqslant a_s^{[k]}$, then $a_{\tilde{t}}^{[k]} \leqslant a_s^{[k]}$ and $\tilde{t} \leqslant s$. Hence,
\[
|b - b_s^{[k]}| \geqslant a_s^{[k]} - |b| = (a_s^{[k]} - a_{\tilde{t}}^{[k]}) + (a_{\tilde{t}}^{[k]} - |b|) \geqslant a_s^{[k]} - a_{\tilde{t}}^{[k]} \geqslant 4\eta_k (s - \tilde{t}).
\]

Summarizing, when $s \notin \{\tilde{t} - 1, \tilde{t}\}$, we deduce the desired inequality \eqref{outsides disc pt}. Lastly, when $s \in \{\tilde{t} - 1, \tilde{t}\}$, i.e., $a_{s-1}^{[k]} \leqslant |b| < a_{s+1}^{[k]}$, we have $|\tilde{t} - s| + 1 \leqslant 2$. In this case, it suffices to show that
\[
|b - b_s^{[k]}| \geqslant \eta_k,
\]
which is obvious since $b \notin \overline{\mathbb{D}}(b_s^{[k]}, \eta_k)$.
\end{proof}

 \subsection{Model Curves}

Our construction is inspired by the strategy presented in \cite{Chen-Huynh-Xie-2023}, with some novel and refined considerations.

\medskip

First, we select a countable and dense subfield of $\mathbb{C}$, for instance, $\overline{\mathbb{Q}}$. We then consider the countable set
\[
\mathscr{R} := \left\{ [p_0(z) : p_1(z) : \dots : p_m(z)] \,:\,
\begin{array}{l}
p_0, p_1, \dots, p_m \in \overline{\mathbb{Q}}[z] \setminus \{0\} \text{ having no common zero in } \mathbb{C}; \\
\deg(p_0) \geqslant \deg(p_i) + 4, \forall i \neq 0.
\end{array}
\right\}
\]
The choice of the number "4" is motivated by analytical reasons that will become clear later (e.g., Proposition~\ref{FHC}).

By employing Runge's approximation theorem, and through the technique of substituting $p_0(z)$ with $p_0(z)\frac{(z+M)^Q}{M^Q}$ for sufficiently large integers $M$ and $Q$, it can be demonstrated that $\mathscr{R}$ is dense in $\mathcal{H}(\mathbb{C}, \mathbb{P}^m)$ with respect to the compact-open topology (see \cite[Lemma~2.1]{Chen-Huynh-Xie-2023}).

We enumerate the elements of $\mathscr{R}$ in any order as follows:
\begin{equation*}
\mathscr{R} = \left\{ \gamma^{[k]}(z) = [p_0^{[k]}(z) : p_1^{[k]}(z) : \dots : p_m^{[k]}(z)] = [1 : \gamma_1^{[k]}(z) : \dots : \gamma_m^{[k]}(z)] \right\}_{k \geqslant 1},
\end{equation*}
where $\gamma_j^{[k]}(z) := p_j^{[k]}(z)/p_0^{[k]}(z)$. Given that $\deg(p_0^{[k]}) \geqslant \deg(p_j^{[k]}) + 4$, there exist sufficiently large radii $R_k > 0$ such that
\begin{equation}\label{decay}
   \left| \gamma_j^{[k]}(z) \right| = \left| \frac{p_j^{[k]}(z)}{p_0^{[k]}(z)} \right|
\leqslant \frac{1}{|z|^3},
\quad \forall |z| \geqslant R_k, \, j=1,\dots,m.
\end{equation}

Let $n_{\gamma^{[k]}}(r, H_0)$ denote the number of zeros of $p_0^{[k]}$ within the disc $\{|z| < r\}$, counted with multiplicities. It is evident that
\begin{equation}\label{num poles}
n_k := \sum_{j=1}^m \#_{\mathsf{multi}} (\gamma_j^{[k]})^{-1}(\infty)
\geqslant n_{\gamma^{[k]}}(r, H_0),
\quad \forall r \geqslant R_k,
\end{equation}
where $\#_{\mathsf{multi}} (\gamma_j^{[k]})^{-1}(\infty)$ represents the number of poles of $\gamma_j^{[k]}$ in $\mathbb{C}$, also counted with multiplicities.

\section{Proof of Theorem~\ref{main-thm-c}}
\label{section 4}

Our desired entire curve will be 
\[ h := [1 : h_1 : \dots : h_m], \]
where (recall \eqref{decay} and \eqref{set center})
\[ h_j := \sum_{k \geqslant 1} h^{[k]}_j(z), \quad h^{[k]}_j := \sum_{s \geqslant 1} \gamma^{[k]}_j(z - b_s^{[k]}). \]
We will show the following properties:
\begin{enumerate}
\item[(1)] (Convergence) $h_j$ is a well-defined meromorphic function for each $j = 1, \dots, m$.
\item[(2)] (Frequent hypercyclicity) $h$ is a common frequently hypercyclic element for all $\mathsf{T}_{e^{\sqrt{-1}\cdot\theta(v)}}$, $v \geqslant 1$.
\item[(3)] (Slow growth rate) There exists some $M > 0$ such that $T_h(r) \leqslant M \cdot r$ for all $r > 0$.
\end{enumerate}

\medskip
First, we prove that $h_j$ is absolutely convergent outside $\bigcup_{s, k \geqslant 1} \overline{\mathbb{D}}(b_s^{[k]}, \eta_k)$.

\begin{Prop}[Convergence]\label{convergence}
For any $b \in \mathbb{C} \setminus \bigcup_{s, k \geqslant 1} \overline{\mathbb{D}}(b_s^{[k]}, \eta_k)$, one has $|h_j(b)| \leqslant 1$.
\end{Prop}

\begin{proof}
For every $k \geqslant 1$, by Proposition~\ref{distance estimate}--(3), we can find some $\tilde{t}_k \in \mathbb{N}$ with $a_{\tilde{t}_k-1}^{[k]} \leqslant |b| < a_{\tilde{t}_k}^{[k]}$ and
\[
|b - b_s^{[k]}| \geqslant \frac{\eta_k}{2} (|\tilde{t}_k - s| + 1) \overset{\eqref{choose eta}}{\geqslant} R_k \quad \text{for all } s \geqslant 1.
\]
Combining this with \eqref{decay}, we get
\begin{equation}\label{control gamma(b)}
\left| \gamma^{[k]}_j(b - b_s^{[k]}) \right| \leqslant \frac{1}{|b - b_s^{[k]}|^3} \leqslant \frac{2^3}{\eta_k^3} \cdot \frac{1}{(|s - \tilde{t}_k| + 1)^3} \quad \text{for all } s \geqslant 1.
\end{equation}
Hence, we can estimate that
\[
\left| h_j^{[k]}(b) \right| \leqslant \sum_{s < \tilde{t}_k} \left| \gamma^{[k]}_j(b - b_s^{[k]}) \right| + \sum_{s \geqslant \tilde{t}_k} \left| \gamma^{[k]}_j(b - b_s^{[k]}) \right| \overset{\eqref{control gamma(b)}}{\leqslant} \frac{2^3}{\eta_k^3} \cdot \sum_{n=1}^{+\infty} \frac{1}{n^3} + \frac{2^3}{\eta_k^3} \cdot \sum_{n=1}^{+\infty} \frac{1}{n^3} \overset{\eqref{choose eta}}{\leqslant} 2^{-k}.
\]
Consequently, we complete the proof using
\[
\left| h_j(b) \right| \leqslant \sum_{k \geqslant 1} \left| h_j^{[k]}(b) \right| \leqslant \sum_{k \geqslant 1} 2^{-k} \leqslant 1.
\]
\end{proof}

Next, we show that $h_j(z)$ is meromorphic on $\bigcup_{s, k \geqslant 1} \overline{\mathbb{D}}(b_s^{[k]}, \eta_k)$, and thereby we verify the frequent hypercyclicity of $h$.

\begin{Prop}[Convergence]\label{FHC}
For any $z \in \overline{\mathbb{D}}(b_s^{[k]}, \Theta^{[k]}(s))$, there holds
\[
\left| h_j(z) - \gamma_j^{[k]}(z - b_s^{[k]}) \right| \leqslant \frac{1}{\Theta^{[k]}(s)}.
\]
\end{Prop}

\begin{proof}
Fix $k, s \geqslant 1$. Recall that $\Theta^{[k]}(s) = l \geqslant \eta_k$, where $|b_s^{[k]}| \in A(k, 2l)$.

For any $t \neq s$, by Proposition~\ref{distance estimate}--(2), we have
\[
\mathsf{dist}(b_t^{[k]}, \overline{\mathbb{D}}(b_s^{[k]}, l)) \geqslant l + \eta_k |t - s| > R_k.
\]
Combining this with \eqref{decay}, we obtain
\[
\left| \gamma^{[k]}_j(z - b_t^{[k]}) \right| \leqslant \frac{1}{(l + \eta_k |t - s|)^3} \quad \text{for all } z \in \overline{\mathbb{D}}(b_s^{[k]}, l).
\]
Thus, we can estimate the difference
\begin{equation}
\label{part 1 estimate}
\left| h_j^{[k]}(z) - \gamma_j^{[k]}(z - b_s^{[k]}) \right| \leqslant \left| \sum_{t > s} \gamma^{[k]}_j(z - b_t^{[k]}) \right| + \left| \sum_{t < s} \gamma^{[k]}_j(z - b_t^{[k]}) \right| \leqslant 2 \sum_{n=1}^{+\infty} \frac{1}{(l + \eta_k n)^3}.
\end{equation}

For any other positive integer $k' \neq k$, by Proposition~\ref{distance estimate}--(2), we can find some integer $\tilde{s} \geqslant 1$ with $a_{\tilde{s}-1}^{[k']} < a_s^{[k]} < a_{\tilde{s}}^{[k']}$, and for any $t \geqslant 1$ there holds
\[
\mathsf{dist}(b_t^{[k']}, \overline{\mathbb{D}}(b_s^{[k]}, l)) \geqslant l + \eta_{k'} (|t - \tilde{s}| + 1) > \eta_{k'} \geqslant R_{k'}.
\]
Applying \eqref{decay}, for any $t \geqslant 1$ and for any $z \in \overline{\mathbb{D}}(b_s^{[k]}, l)$, there holds
\[
\left| \gamma^{[k']}_j(z - b_t^{[k']}) \right| \leqslant \frac{1}{(l + \eta_{k'} (|t - \tilde{s}| + 1))^3},
\]
which implies
\[
\left| h_j^{[k']}(z) \right| \leqslant \left| \sum_{t \geqslant \tilde{s}} \gamma^{[k']}_j(z - b_t^{[k']}) \right| + \left| \sum_{t < \tilde{s}} \gamma^{[k']}_j(z - b_t^{[k']}) \right| \leqslant 2 \sum_{n=1}^{+\infty} \frac{1}{(l + \eta_{k'} n)^3}.
\]
This inequality, combined with \eqref{part 1 estimate}, gives
\[
\left| h_j(z) - \gamma_j^{[k]}(z - b_s^{[k]}) \right| \leqslant \left| h_j^{[k]}(z) - \gamma_j^{[k]}(z - b_s^{[k]}) \right| + \sum_{k' \neq k} \left| h_j^{[k']}(z) \right| \leqslant \sum_{a=1}^{+\infty} \sum_{n=1}^{+\infty} \frac{2}{(l + \eta_a n)^3}.
\]
Now, 
\[
\sum_{n=1}^{+\infty} \frac{2}{(l + \eta_a n)^3} \leqslant \frac{2}{(l + \eta_a)^3} + \sum_{n=2}^{+\infty} \int_{x=n-1}^{n} \frac{2}{(l + \eta_a x)^3} \, \mathrm{d}x \leqslant \frac{2}{\eta_a (l + \eta_a)^2} + \int_{x=1}^{+\infty} \frac{2}{(l + \eta_a x)^3} \, \mathrm{d}x = \frac{3}{\eta_a (l + \eta_a)^2}.
\]
Thus the proof
follows from the following estimate:
\[
\sum_{a=1}^{+\infty} \sum_{n=1}^{+\infty} \frac{2}{(l + \eta_a n)^3} \leqslant \sum_{a=1}^{+\infty} \frac{3}{\eta_a (l + \eta_a)^2} \leqslant \sum_{a=1}^{+\infty} \frac{3}{\eta_a^2} \cdot \frac{1}{l} \overset{\eqref{choose eta}}{\leqslant} \sum_{a=1}^{+\infty} \frac{3}{9^a} \cdot \frac{1}{l} \leqslant \frac{1}{l}.
\]
\end{proof}

\begin{Prop}[Frequent Hypercyclicity]
This entire curve $h \in \bigcap_{v \geqslant 1} \mathsf{FHC}(\mathsf{T}_{e^{\sqrt{-1}\cdot\theta(v)}})$.
\end{Prop}

\proof
By Proposition~\ref{FHC}, we have
\begin{equation}\label{FHC2}
\sup_{|z|\leqslant \Theta^{[k]}(s)}\left|h_j(z+b_s^{[k]})-\gamma_j^{[k]}(z)\right|
\leqslant 
\frac{1}{\Theta^{[k]}(s)},
\quad \forall j = 1, \dots, m.    
\end{equation}
For each element $\theta(v)$ in $E$, we have
\[
\{b^{[k]}_s\}_{s\geqslant 1}\cap \mathbb{N}\cdot e^{\sqrt{-1}\cdot\theta(v)}
\overset{\eqref{center2}}{=}
\mathcal{B}^{[k]}_v=\mathcal{A}^{[k]}_v\cdot e^{\sqrt{-1}\cdot\theta(v)}
\overset{\eqref{center1}}{=}
\bigcup_
{
l\geqslant \eta_k,\, l\in\varphi_1^{-1}(v)
} A(k,2l)\cdot e^{\sqrt{-1}\cdot\theta(v)}.
\]
The inequality~\eqref{FHC2} implies that for any $n \in A(k,2l)$, there holds
\[
\sup_{|z|\leqslant l} \left|\mathsf{T}^n_{e^{\sqrt{-1}\cdot\theta(v)}}\big(h_j(z)\big)-\gamma_j^{[k]}(z)\right|\leqslant \frac{1}{l},
\quad \forall j = 1, \dots, m.
\]
Noting that 
$\{\gamma^{[k]}\}_{k\geqslant 1}$ is dense in $ \mathcal{H}(\mathbb{C},\mathbb{CP}^m)$, 
and that every
$A(k,2l)$ has positive lower density, we finish the proof. 
\qed

\medskip

\begin{Prop}[Slow Growth Rate]\label{slow growth}
The entire curve $h$ grows slowly in the sense that
\[
T_h(r) \leqslant O(r), \quad \forall r > 1.
\]
\end{Prop}

Here and from now on, we write a term as $O(r)$ when it is less than $C \cdot r$
for some uniform constant $C > 0$ and for all $r$ according to the context.

\begin{proof}
There are two cases to be treated, depending on the value of $r > 1$.

\noindent{\bf Case 1:} The circle $\partial \mathbb{D}_r := \{|z| = r\}$ does not intersect any open disc in $\{{\mathbb{D}}(b_s^{[k]},\eta_k)\}_{s\geqslant 1, k\geqslant 1}$. In this case, Proposition~\ref{convergence} guarantees that, for any $z \in \partial \mathbb{D}_r$ and $1 \leqslant j \leqslant m$, we have $|h_j(z)| \leqslant 1$. Hence the proximity function
\begin{equation}\label{mhr}
m_h(r,H_0) \leqslant \log\sqrt{m+1} = O(1).
\end{equation}

For the counting function, we have

\begin{Cla}\label{counting}
$
n_h(t,H_0) \leqslant O(t)
$ for all $t \geqslant 1$.
\end{Cla}

\smallskip
Assuming the above claim at the moment, we conclude that
\begin{equation}
\label{Nhr}
N_h(r,H_0) = \int_{t=1}^r \frac{n_h(t,H_0)}{t} \, \text{d}t + O(1) \leqslant O(r) + O(1).
\end{equation}
Summing up \eqref{mhr} and \eqref{Nhr}, and using \eqref{fmt}, we get
\begin{align}\label{ak+rk}
T_h(r) = N_h(r,H_0) + m_h(r,H_0) + O(1) &\leqslant O(r).
\end{align}

\noindent{\bf Case 2:} The circle $\{|z| = r\}$ intersects some disc ${\mathbb{D}}(b_s^{[k]},\eta_k)$. In this case, $|b_s^{[k]}| - \eta_k < r < |b_s^{[k]}| + \eta_k$. From \eqref{farfrom0}, we obtain $\eta_k \leqslant |b_s^{[k]}| - \eta_k$ and
\begin{equation}\label{enlarge radius}
r < |b_s^{[k]}| + \eta_k = (|b_s^{[k]}| - \eta_k) + 2\eta_k \leqslant 3(|b_s^{[k]}| - \eta_k) \leqslant 3r.    
\end{equation}
By Proposition~\ref{distance estimate}--(1), the circle $\{|z| = |b_s^{[k]}| + \eta_k\}$ does not intersect any open disc $\{{\mathbb{D}}(b_s^{[k]},\eta_k)\}_{s\geqslant 1, k\geqslant 1}$.
Since $T_{h}(\cdot)$ is an increasing function, we conclude that
\[
T_{h}(r) \overset{\eqref{enlarge radius}}{\leqslant} T_{h}(|b_s^{[k]}| + \eta_k) \overset{\eqref{ak+rk}}{\leqslant} O(|b_s^{[k]}| + \eta_k) + O(1) \overset{\eqref{enlarge radius}}{\leqslant} O(r).
\]
This finishes the proof.
\end{proof}

\begin{proof}[Proof of the Claim~\ref{counting}]
By Proposition~\ref{convergence}, $h_j(z)$ has no poles outside the discs $\{\mathbb{D}(b_s^{[k]},\eta_k)\}_{s,k\geqslant 1}$, and the number of poles in each $\mathbb{D}(b_s^{[k]},\eta_k)$ is (recall $\eta_k > R_k$)
\[
\#_{\mathsf{multi}} h_j^{-1}(\infty) \cap \mathbb{D}(b_s^{[k]},\eta_k) = \#_{\mathsf{multi}} (\gamma_j^{[k]})^{-1}(\infty) \cap \mathbb{D}(0,\eta_k) \overset{\eqref{decay}}{=} \#_{\mathsf{multi}} (\gamma_j^{[k]})^{-1}(\infty).
\]
Set $\cN_t := \{b_s^{[k]} : s, k \geqslant 1, \mathbb{D}(b_s^{[k]},\eta_k) \cap \mathbb{D}(0,t) \neq \emptyset\}$. Then
\begin{equation}
\label{count zeros part}
n_h(t,H_0) \leqslant \sum_{b_s^{[k]} \in \cN_t} \sum_{j=1}^m \#_{\mathsf{multi}} h_j^{-1}(\infty) \cap \mathbb{D}(b_s^{[k]},\eta_k) = \sum_{b_s^{[k]} \in \cN_t} \sum_{j=1}^m \#_{\mathsf{multi}} (\gamma_j^{[k]})^{-1}(\infty) \overset{\eqref{num poles}}{=} \sum_{b_s^{[k]} \in \cN_t} n_k.
\end{equation}

Note that by our construction
\begin{equation}\label{planar obs}
\bigcup_{b_s^{[k]} \in \cN_t} \{|b^{[k]}_s| - \eta_k < |z| < |b^{[k]}_s| + \eta_k\} \subset \mathbb{D}(0,3t).    
\end{equation}
Indeed, for any $b_s^{[k]} \in \cN_t$,
we have $t > |b_s^{[k]}| - \eta_k \geqslant \eta_k$ by \eqref{farfrom0}.
Hence $|b_s^{[k]}| + \eta_k = (|b_s^{[k]}| - \eta_k) + 2\eta_k < 3t$, and \eqref{planar obs} follows.

By Proposition~\ref{distance estimate}--(1),
the annuli on the left-hand-side of \eqref{planar obs} are pairwise disjoint. Hence
\[
3t \geqslant \sum_{b_s^{[k]} \in \mathcal{N}_t} 2\eta_k \overset{\eqref{choose eta}}{\geqslant} \sum_{b_s^{[k]} \in \mathcal{N}_t} 2 n_k.
\]
Therefore $n_h(t,H_0) \overset{\eqref{count zeros part}}{\leqslant} \sum_{b_s^{[k]} \in \mathcal{N}_t} n_k \leqslant O(t)$.
\end{proof}

\begin{proof}[Proof of Theorem~\ref{main-thm-c}]
Since
$
T_h(r) = O(r^2) 
$ 
for very small $r > 0$, we have
$\lim_{r \to 0^+} T_h(r)/r = 0$.
Noting that $T_h(r)/r$ is a continuous function with respect to $r > 0$,
and that $\sup_{r > 1} T_h(r)/r < +\infty$
by Proposition~\ref{slow growth}, it is clear that  
$M := \sup_{r > 0} T_h(r)/r < +\infty$. Whence
\begin{equation}
\label{T(r)<Mr}
T_h(r) \leqslant M \cdot r,
\quad
\forall r > 0.
\end{equation}

Given $\epsilon > 0$, we now define a new entire curve
\begin{equation}
\label{magic trick}
\tilde{h}(\bullet) := h\left(\frac{\epsilon}{M} \cdot \bullet\right): \mathbb{C} \rightarrow \mathbb{P}^m,
\end{equation}
which clearly lies in $\bigcap_{v \geqslant 1} \mathsf{FHC}(\mathsf{T}_{\frac{M}{\epsilon} \cdot e^{\sqrt{-1}\cdot\theta(v)}})$. 
It is well-known that $\mathsf{FHC}(\mathsf{T}_{e^{\sqrt{-1}\cdot\theta}}) = \mathsf{FHC}(\mathsf{T}_{R \cdot e^{\sqrt{-1}\cdot\theta}})$ for any $R > 0$ (cf. e.g.~\cite[Lemma~4.1]{Guo-Xie-2}). Therefore $\tilde{h} \in \bigcap_{v \geqslant 1, R > 0} \mathsf{FHC}(\mathsf{T}_{R \cdot e^{\sqrt{-1}\cdot\theta(v)}})$.

Lastly, we estimate the growth rate of $\tilde{h}$ by direct calculation
\begin{align*}
T_{\tilde{h}}(r) &= \int_{t=0}^r \frac{\mathrm{d}t}{t} \int_{|z|<t} \tilde{h}^*\omega_{\sf FS} = \int_{t=0}^{\frac{\epsilon}{M} \cdot r} \frac{\mathrm{d}t}{t} \int_{|z|<t } {h}^*\omega_{\sf FS} \\
&= T_h\left(\frac{\epsilon}{M} \cdot r\right) \overset{\eqref{T(r)<Mr}}{\leqslant} M \cdot \frac{\epsilon}{M} \cdot r = \epsilon \cdot r,
\quad
\forall r > 0.    
\end{align*}
Hence we are done.
\end{proof}

\bigskip

\noindent\textbf{Acknowledgements}

\smallskip\noindent
We thank Dinh Tuan Huynh for inspiring conversations and comments. We thank Yi C. Huang for careful reading of the manuscript. We thank the anonymous referee for their careful reading and constructive suggestions, which have enhanced the clarity of the manuscript.

\bigskip

\noindent\textbf{Funding}

\smallskip\noindent
Song-Yan Xie was partially supported by: the National Key R\&D Program of China under Grants No. 2021YFA1003100 and No. 2023YFA1010500; the National Natural Science Foundation of China under Grants No. 12288201 and No. 12471081; and the Xiaomi Young Talents Program.
Zhangchi Chen was partially supported by: the National Natural Science Foundation of China (Grant No. 12501104); the Science and Technology Commission of Shanghai Municipality (Grant No. 22DZ2229014); the Shanghai Sailing Program (Grant No. 24YF2709900); the Shanghai Pujiang Program (Grant No. 24PJA023); the Labex CEMPI (ANR-11-LABX-0007-01); and the project QuaSiDy (ANR-21-CE40-0016).

\bigskip


\setlength\parindent{0em}

{\scriptsize Zhangchi Chen
\\
School of Mathematical Sciences, Key Laboratory of MEA (Ministry of Education) and Shanghai Key Laboratory of PMMP, East China Normal University, Shanghai 200241, China
}

{\bf\scriptsize zcchen@math.ecnu.edu.cn}

\medskip

{\scriptsize Bin Guo
\\
Academy of Mathematics and Systems Sciences, Chinese Academy of Sciences, Beijing 100190, China}
\\
{\bf\scriptsize guobin181@mails.ucas.ac.cn}

\medskip

{\scriptsize Song-Yan Xie
\\
State Key Laboratory of Mathematical Sciences, Academy of Mathematics and Systems Science, Chinese Academy of Sciences, Beijing 100190, China;  School of Mathematical Sciences, University of Chinese Academy of Sciences, Beijing 100049, China.}
\\
{\bf\scriptsize xiesongyan@amss.ac.cn}


\begin{thebibliography}{XL}{\scriptsize

\medskip

{\bf\bibitem{BCP-2021}
{\rm Bayart}} F.; {\rm Costa Júnior} F.; {\rm Papathanasiou} D.:
{\em Baire theorem and hypercyclic algebras.}
Adv. Math. {\bf 376} (2021), Paper No. 107419, 58 pp.

\medskip

{\bf\bibitem{Bayart-Grivaux-2006}
{\rm Bayart}} F.; {\rm Grivaux} S.:
{\em Frequently hypercyclic operators.}
Trans. Amer. Math. Soc., {\bf 358} (2006), no. 11, 5083–5117.

\medskip

{\bf\bibitem{BGMM-2024}
{\rm Bayart}} F.; {\rm Grivaux} S.; {\rm Matheron} E.; {\rm Menet} Q.: 
{\em Hereditarily frequently hypercyclic operators and disjoint frequent hypercyclicity}, arXiv:2409.07103.

\medskip

{\bf\bibitem{Bayart-Matheron-2009}
{\rm Bayart}} F.; {\rm Matheron} É.:
{\em Dynamics of linear operators.}
Cambridge Tracts in Mathematics, {\bf 179}. Cambridge University Press, Cambridge, (2009). xiv+337 pp. ISBN: 978-0-521-51496-5

\medskip


{\bf\bibitem{Birkhoff-1929}
{\rm Birkhoff}} G. D.:
{\em D\'emonstration d'un th\'eor\`eme \'el\'ementaire sur les fonctions enti\`eres.}
C.R. Acad. Sci. Paris, {\bf 189} (1929), 473--475.

\medskip

{\bf\bibitem{BBGE}
{\rm Blasco}} O; {\rm Bonilla} A.; 
{\rm  Grosse-Erdmann} K.-G.:
{\em Rate of growth of frequently hypercyclic functions.} 
Proc. Edinb. Math. Soc. (2), {\bf 53} (2010) no.1, 39--59.

\medskip

{\bf\bibitem{Chen-Huynh-Xie-2023}
{\rm Chen}} Z.; {\rm Huynh} D. T.; {\rm Xie} S.-Y.:
{\em Universal entire curves in projective spaces with slow growth.}
J. Geom. Anal., {\bf 33} 308 (2023).

\medskip

{\bf\bibitem{Dinh-Sibony-2020}
{\rm Dinh}} T-C.; {\rm Sibony} N.;
{\em Some open problems on holomorphic foliation theory.}
Acta Math. Vietnam., {\bf 45} (2020), no.1, 103--112.


\medskip

{\bf\bibitem{EM-2019}
{\rm Ernst}} R.; {\rm Mouze} A.;
{\em A quantitative interpretation of the frequent hypercyclicity criterion.}
Ergodic Theory and Dynamical Systems.,   {\bf 39 (4)} (2019), 898--924.

\medskip

{\bf\bibitem{EM-2021}
{\rm Ernst}} R.; {\rm Mouze} A.;
{\em Frequent universality criterion and densities.}
Ergodic Theory and Dynamical Systems.,   {\bf 41(3)} (2021), 846--868. 

\medskip

{\bf\bibitem{Forstneric-2023}
{\rm Forstneri\v{c}}}, F.:
{\em Recent developments on Oka manifolds.}
Indag. Math. {\bf 2}, (2023), vol 34, 367--417.



\medskip

{\bf\bibitem{Grosse-Erdmann-Manguillot-2011}
{\rm Grosse-Erdmann}} K.; {\rm Manguillot} A.P.:
{\em Linear chaos.}
Springer London; Universitext; 2011.

\medskip

{\bf\bibitem{Guo-Xie}
{\rm Guo}} B.; 
{\rm Xie} S.-Y.:
{\em Universal holomorphic maps with slow growth I: an algorithm.}
Math. Ann., {\bf 389} (2024), 3349--3378.

\medskip

{\bf\bibitem{Guo-Xie-2}
{\rm Guo}} B.; 
{\rm Xie} S.-Y.:
{\em Universal holomorphic maps with slow growth II:  functional analysis methods}, arXiv:2310.06561, to appear in Indiana Univ. Math. J.

\medskip

{\bf\bibitem{Noguchi-Winkelmann-2014}
{\rm Noguchi}} J.; {\rm Winkelmann} J.:
{\em Nevanlinna theory in several complex variables and Diophantine approximation.}
Grundlehren der mathematischen Wissenschaften 350. Springer, Tokyo, (2014). xiv+416 pp.

\medskip

{\bf\bibitem{Ru-2021}
{\rm Ru}} M.:
{\em Nevanlinna theory and its relation to Diophantine approximation.}
Second edition [of MR1850002]. World Scientific Publishing Co. Pte. Ltd., Hackensack, NJ, (2021). xvi+426 pp. 

\medskip


{\bf\bibitem{Voronin}
{\rm Voronin}}, S. M.:
 {\em Theorem on the ``Universality'' of the Riemann Zeta Function.} Izv. Akad. Nauk SSSR, Ser. Mat., 39, no. 3, 475--486, 1975.

}
\end{thebibliography}
\end{document}